# NUMERICAL SOLUTION OF THE SEVENTH ORDER BOUNDARY VALUE PROBLEMS USING B-SPLINE METHOD

## A PREPRINT

**Maryam Khazaei**
Department of Computer Science
University of California, Merced
mkhazaeipool@ucmerced.edu

**Yeganeh Karamipour**
Department of Mathematics
Iran University of Science and Technology
y.karamipour@gmail.com

September 13, 2021

## ABSTRACT

We develop a numerical method for solving the boundary value problem of The Linear Seventh Ordinary Boundary Value Problem by using seventh degree B-Spline function. Formulation is based on particular terms of order of seventh order boundary value problem. We obtain Septic B-Spline formulation and the Collocation B-spline method is formulated as an approximation solution. We apply the presented method to solve an example of seventh-order boundary value problem which the result show that there is an agreement between approximate solutions and exact solutions. Resulting low absolute errors show that the presented numerical method is effective for solving high order boundary value problems. Finally, a general conclusion has been included.

*Keywords* Spline and B-spline functions · Seventh order Boundary value problem · Septic B-spline collocation method

## 1 Introduction

We consider the linear seventh order differential equation

$$y^7(x) = g(x)y(x) + q(x) \tag{1}$$

with the boundary conditions

$$y(a) = k_1, \quad y'(a) = k_2, \quad y''(a) = k_3, \quad y'''(a) = k_4, \quad y(b) = k_5, \quad y'(b) = k_6, \quad y''(b) = k_7, \tag{2}$$

Obtaining the analytical solution of this problem is difficult in general. Scientists have used spline functions more commonly due to their useful application in applied mathematics and engineering. Many fields including Computational Mathematics, Mathematical Physics, and Mechanics have played an important role in studying the application of B-Spline (8) methods.

Higher order differential equations have been solved by many numerical and semi-analytical methods. The solutions of differential equations based on spline functions contain high accuracy. There are comprehensive study with scientists to discuss the numerical solution of linear and nonlinear boundary value problems. For instance, quadratic, cubic, quartic, quintic, sextic, septic and higher degree spline techniques have been used for this purpose. Kumar and Srivastava (9) have given a survey research on cubic, quintic and sextic polynomial and nonpolynomial spline techniques for solving boundary value problems.

For finding the numerical solution of a particular case of seventh order linear boundary value problems where some terms of the boundary value problem are zero, we use B-Spline Basis and collocation method which is considered as approximation solution. Considering a generic case with a general form of Boundary value problem including term of B-spline basis functions with higher derivatives in the general form of the boundary value problem results a dense linear system which takes more computational costs with each added terms. Thus, resulting a sparse and banded matrix in the presented method has influenced to have an efficient method in terms of time and cost. Comparing with some higher order B-spline techniques, the presented method shows efficient numerical results with less absolute errors. The method could also be extended to solve higher order boundary value problems.

This paper is organized as follows: In Sections 2 and 3, an introduction and the fundamental definitions of spline, B-spline functions and obtained B-spline basis functions are presented. In Section 3, we consider collocation method based on B-spline for the linear seventh order boundary value problem. In Section 4, we have discussed about the numerical results. In Section 5, the conclusion and further developments are given.

## 2 Definitions

### 2.1 Spline

A spline is described as a piecewise polynomial function that is being as smooth as it can be without reducing it to a polynomial (1). In general, Spline is a piecewise polynomial function defined in the interval $[a, b]$, such that there exists a decomposition of $D$ into sub-regions in each of which the function is a polynomial of some degree $k$. The term "spline" is used to refer to a wide class of functions that are used in applications requiring data interpolation or smoothing. A function $S(x)$ is a spline of degree $k$ on $[a, b]$ if

$$S \in C^{k-1}[a, b]$$

$$a = t_0 < t_1 < ... < t_n = b$$

and

$$S(x) = \begin{cases} S_0(x), & t_0 < x < t_1, \\ S_1(x), & t_1 < x < t_2, \\ S_{n-1}(x), & t_{n-1} < x < t_n \end{cases} \tag{3}$$

where

$$S_i(x) \in P^k, \quad i = 0, 1, ..., n-1$$

### 2.2 B-Spline

The B-spline of degree $k$ is denoted by $\varphi_i^k(x)$, where $i \in Z$,

We define one more parameter to define B-spline basis functions. The $i - th$ B-spline basis function of degree p, written as $N_{i,p(u)}$, is defined recursively as follows:

$$N_{i,p(u)} = \begin{cases} 1, & \text{if } u_{i+1} > u \geq u_i \\ 0, & \text{otherwise} \end{cases} \tag{4}$$

$$N_{i,p(u)} = \frac{u - u_i}{u_{i+p} - u_i} N_{i,p-1}(u) + \frac{u_{i+p+1} - u}{u_{i+p+1} - u_{i+1}} N_{i+1,p-1}(u) \tag{5}$$

B-spline, or basis spline, is a Spline function that has minimal support with respect to a given degree, smoothness, and domain partition. Any Spline function of a given degree can be expressed as a linear combination of B-splines of that degree. Cardinal B-Splines have knots that are equidistant from each other. Here, we denote the B-spline of degree $k$ by $B_i^k(x)$, where $i$ is an element in $\mathbb{Z}$ with the following properties:

1. $N_{i,p}(u)$ is a degree p polynomial in $u$.

2. Non-negativity: for all $i, p$ and $u$, $N_{i,p}(u)$ is non-negative.

3. Local support: $N_{i,p}(u)$ is a nonzero polynomial on $[ui, u_{i+p+1})$

4. At most $p + 1$ degree $p$ basis functions are nonzero on any span $[ui, u_{i+p+1}]$, namely: $N_{i-p,p}(u), N_{i-p+1,p}(u), N_{i-p+2,p}(u), ..., N_{i,p}(u)$. This property shows that the following basis functions are nonzero on $[ui, u_{i+p+1})$.

$$N_{i-p,p}(u), \quad N_{i-p+1,p}(u), \quad N_{i-p+2,p}(u), \quad ..., N_{i,p}(u)$$

5. Partition of unity: The sum of all nonzero degree p basis functions on span $[u_i, u_{i+p+1})$ is 1 which states that the sum of these $p+1$ basis functions is 1.

**Obtaining the Septic B-Spline**

We consider equally-spaced knots of a partition $\pi : a = x_0 < x_1 < ... < x_n$ on $[a, b]$. The alternative approach for deriving the B-splines which are more applicable with respect to the recurrence relation for the formulations of B-splines of higher degrees. At first, we recall that the $kth$ forward difference $f(x_0)$ of a given function $f(x)$ at $x_0$, which is defined recursively by the following (10) and (11):

$$\triangle f(x_0) = \triangle f(x_1) - \triangle f(x_0), \quad \triangle^{k+1} f(x_0) = \triangle^k f(x_1) - \triangle^k f(x_0) \tag{6}$$

**Definition:** The function $(x - t)m_+$

$$(x - t)^m_+ = \begin{cases} (x - t)^m, & x \leq t_0 \\ 0, & x < t \end{cases} \tag{7}$$

It is clear that $(x - t)^m_+$ is $(m - 1)$ times continuously-differentiable with respect to $t$ and $x$. The B-spline of order m is defined as follows:

$$B_i^m(t) = \frac{1}{h^m} \sum_{j=0}^{m+1} (-1)^{m+1-j} (x_{i-2+j} - t)^m_+ = \frac{1}{h^m} \triangle^{m+1} (x_{i-2} - t)^m_+ \tag{8}$$

Considering various values of $m$. Let $m = 1$, turns out the B-spline of various orders.

$$\frac{1}{h^1} \triangle^2 (x_{i-2} - t)^1_+ = \frac{1}{h^1} [(x_{i-2} - t)_+ - 2(x_{i-1} - t)_+ + (x_i - t)_+]$$

$$B_i^m(t) = \begin{cases} (x_{i-2} - t) - 2(x_{i-1} - t), & x_{i-2} < t \leq x_{i-1} \\ (x_i - t), & x_{i-1} < t \leq x_i \\ 0, & otherwise \end{cases} \tag{9}$$

By considering different value for $m$, we can get different degree of B-Spline. To obtain Septic B-spline basis functions, we need to partition the interval $[a, b]$ while we need to choose mesh points such that $a = x_0, .... b = x_n$ and $x_i = a + ih, i = 0, 1, ...., N$ where $h = (b - a)/n$.

$$B_i(x) = \frac{1}{h^7} \begin{cases} (X - X_{j-4})^7, & X \in [X_{j-4}, X_{j-3}] \\ (X - X_{j-4})^7 - 8(X - X_{j-3})^7, & X \in [X_{j-3}, X_{j-2}] \\ (X - X_{j-4})^7 - 8(X - X_{j-3})^7 + 28(X - X_{j-2})^7, & X \in [X_{j-2}, X_{j-1}] \\ (X - X_{j-4})^7 - 8(X - X_{j-3})^7 + 28(X - X_{j-2})^7 - 56(X - X_{j-1})^7, & X \in [X_{j-1}, X_j] \\ (X_{j+4} - X)^7 - 8(Xj + 3 - X)^7 + 28(X_{j+2} - X)^7 - 56(X_{j+1} - X)^7, & X \in [X_j, X_{j+1}] \\ (X_{j+4} - X)^7 - 8(Xj + 3 - X)^7 + 28(X_{j+2} - X)^7, & X \in [X_{j+1}, X_{j+2}] \\ (X_{j+4} - X)^7 - 8(Xj + 3 - X)^7, & X \in [X_{j+2}, X_{j+3}] \\ (X_{j+4} - X)^7, & X \in [X_{j+2}, X_{j+3}] \\ 0, \end{cases}$$

let $S_7[\pi]$ be the space of continuously-differentiable, piecewise seventh degree polynomials on $\pi$. Using (3), the B-splines and its derivatives are defined in the table $1.a, 1.b, 1.c$ respectively. Table 1 with the value of $B_i^{''}(x)$ and $B_i^4(x)$, and $Table 2$ with the value of $B_i^5(x)$ and $B_i^6(x)$ and finally, $Table 3$ The value of $B_i^7(x)$ and $B_i^6(x)$ are shown respectively. We also point out that $B_i(x) = 0$ for $x < x_{i-4}$ and $x > x_{i+4}$.

Table 1: Values of $B_i^3(x)$ and $B_i^4(x)$

| $x_{i-4}$ | $x_{i-3}$ | $x_{i-2}$ | $x_{i-1}$ | $x_i$ | $x_{i+1}$ | $x_{i+2}$ | $x_{i+3}$ | $x_{i+4}$ |
|---|---|---|---|---|---|---|---|---|
| 0 | $210/h^3$ | $1680/h^3$ | $-3990/h^3$ | 0 | $3990/h^3$ | $-1680/h^3$ | $-210/h^3$ | 0 |
| 0 | $840/h^4$ | 0 | $-7560/h^4$ | $13440/h^4$ | $-7560/h^4$ | 0 | $840/h^4$ | 0 |

Table 2: Values of $B_i^5(x)$ and $B_i^6(x)$

| $x_{i-4}$ | $x_{i-3}$ | $x_{i-2}$ | $x_{i-1}$ | $x_i$ | $x_{i+1}$ | $x_{i+2}$ | $x_{i+3}$ | $x_{i+4}$ |
|---|---|---|---|---|---|---|---|---|
| 0 | $2520/h^5$ | $-10080/h^5$ | $12600/h^5$ | 0 | $-12600/h^5$ | $10080/h^5$ | $-2520/h^5$ | 0 |
| 0 | $5040/h^6$ | $-30240/h^6$ | $75600/h^6$ | $-100800/h^6$ | $75600/h^6$ | $-30240/h^6$ | $5040/h^6$ | 0 |

## 3 Collocation method based on B-spline for linear seven order boundary value problem

The approximate solution of equation (1.1) by Collocation method based on septic B-spline is as follows

$$y(x) = \sum_{j=-7}^{n} \alpha_j \beta_j(x) \tag{10}$$

Where $\alpha_j$, are unknown real coefficients to be determine and $B_j(x)$ are seven degree B-spline function and let $x_0, x_1, ..., x_n$ be $n + 1$ grid points in the interval $[a, b]$, so that

$$x_i = a + ih, i = 0, \quad 1, ...n, \quad x_0 = a, \quad x_n = b, \quad h = (b-a)/n$$

with substituting (12) in (1), we get the following

$$[\sum_{j=-7}^{n} \alpha_j (\beta_j^7(x_i) - g(x_i)\beta_j(x_i)) - q(x_i)] = 0, \quad i = 0, 1, .., n \tag{11}$$

We use of boundary conditions (2) and (12) to make the following system to evaluate the $\alpha_j$ so that we get the following:

$$[\sum_{j=-7}^{n} \alpha_j (\beta_j^7(x_i) - g(x_i)\beta_j(x_i)) - q(x_i)] = 0, \quad i = 0, 1, .., n \tag{12}$$

$$\sum_{j=-7}^{n} \alpha_j \beta_j(x_i) = y(a), \qquad i = n+1 \tag{13}$$

$$\sum_{j=-7}^{n} \alpha_j \beta_j'(x_i) = y'(a), \qquad i = n+2 \tag{14}$$

$$\sum_{j=-7}^{n} \alpha_j \beta_j''(x_i) = y''(a), \qquad i = n+3 \tag{15}$$

Table 3: Values of $B_i^7(x)$

| $x_{i-4}$ | $x_{i-3}$ | $x_{i-2}$ | $x_{i-1}$ | $x_i$ | $x_{i+1}$ | $x_{i+2}$ | $x_{i+3}$ | $x_{i+4}$ |
|---|---|---|---|---|---|---|---|---|
| 0 | $5040/h^7$ | $-35280/h^7$ | $105840/h^7$ | $-176400/h^7$ | $-176400/h^7$ | $105840/h^7$ | $-35280/h^7$ | 0 |

Table 4: The values of Exact solution, Approximate Solution and Max Absolute Errors for problem 1

| n | $Y_i$ | $y_i$ | Max absolute Error |
|----|-----------|--------------|-------------------|
| 20 | 0.35094702 | 0.33068119 | 0.020265825 |
| 40 | 0.80970360 | 0.80296795 | 0.00673565 |
| 60 | 0.99985955 | 0.9921557165 | 0.00602099 |
| 80 | 1.00000000 | 0.99969918 | 0.00030082 |

$$\sum_{j=-7}^{n} \alpha_j \beta_j'''(x_i) = y'''(a), \quad i = n+4 \tag{16}$$

$$\sum_{j=-7}^{n} \alpha_j \beta_j''(x_i) = y''(b), \quad i = n+5 \tag{17}$$

$$\sum_{j=-7}^{n} \alpha_j \beta_j'(x_i) = y'(b), \quad i = n+6 \tag{18}$$

$$\sum_{j=-7}^{n} \alpha_j \beta_j(x_i) = y(b), \quad i = n+7 \tag{19}$$

Using Table $(1.a - c)$ and $(3.4) to (3.11)$, we evaluate the values of B-Spline functions at the knots $x_i, i = 0, .., n$ A system $A$ of $n + 8$ linear equations in the $(n+8)$ unknowns $\alpha_{-7}, \alpha_{-6}, ...., \alpha_n$ is obtained. This system can be written in a matrix-vector form as follows

$$A\alpha = X \tag{20}$$

where

$$\alpha = [\alpha_{-7}, \alpha_{-6}, ...., \alpha_n]^T$$

and

$$X = [q(x_0), ..., q(x_n), y(a), y'(a), y''(a), y'''(a), y''(b), y'(b), y(b)]^T$$

The matrix $A$ is a linear matrix which is sparse and considered as the banded matrix which is a sparse matrix whose nonzero entries are confined to a diagonal band, comprising the main diagonal and zero or more diagonals on either side. These properties of matrix $A$ leads more efficiency and less computational cost for the method as well as the algorithm of Cholesky or the $LDL^T$ factorization can take advantage of the structure banded matrices. By solving the linear system, the values of

$$\alpha = [\alpha_{-7}, \alpha_{-6}, ...., \alpha_n]^T$$

can be obtained, then the approximate solution of equation (3.1) will be evaluated by using (3.12).

## 4 Numerical results

In this section, we solve a linear problems that will be solved by septic B-Spline functions.

Example: First, We consider the linear boundary value problem (4).

$$y^7(x) = y(x) - 7e^x, \quad 0 \le x \le 1 \tag{21}$$

with the following boundary conditions

$$y(0) = 1, \quad y(1) = 0, \quad y'(1) = e, \quad y''(1) = 2e, \quad y'(0) = 0, \quad y''(0) = 1, \quad y'''(0) = 2, \tag{22}$$

The domain $[0, 1]$ is divided into $20, 40, 60, 80$ equal sub-intervals. we define mesh-points $x_i$ as $x_i = x_0 + ih$ where $h$ is $h = (b = x_n - a = x_0)/n$. Maximum absolute error and approximate solution $y_i$ coming from the presented method and exact solution $Y_i$ for this problem are shown in the table 4, which the analytical solution is $(1-x)e^x$. A comparison for the absolute error of (7) and ninth degree spline approximation for Example 1 at $h = 0.1$, turns out that (7) has the larger absolute error rather than our presented method.

## 5   Conclusion

For finding the numerical solution of a particular case of seventh order linear boundary value problems where some terms of the BVP are zero, we use B-Spline Basis and Collocation method which is considered as an approximation solution. We evaluate the solutions for problem by the seventh degree B-spline. We end up with a banded sparse matrix with the property of being diagonal that takes less computational cost than some other methods (6) which they consider general form of Boundary value problem. The reason is that there will be more terms of B-spline basis functions with higher derivatives in the general form of the boundary value problem. Then, the linear system form will be contained a more dense matrix which takes more computational cost with each added terms. Therefore, the presented method will be required less running time. Numerical results show the capability and efficiency of the present method. Numerical results showed that the method achieved forth-order accuracy. One of the importance of this method is that we can approximate the solution at every point of the range of domain. Numerical results obtained by the present method are in good agreement with the exact solutions or numerical solutions available in the literature. We can extend the degree of the proposed method to solve the higher (i.e., more than 7th) order boundary value problem. The presented method also could be used for nonlinear form of seventh order Boundary value problems with the same B-Spline basis functions.